\documentclass[11pt]{article}
\usepackage{amsthm,enumerate}
\renewcommand{\baselinestretch}{1.2}
\newcommand{\dated}{\mbox{} \hfill {\small [{\tt \today}]}} \usepackage{amsmath,amssymb,amsfonts,diagrams}
\newenvironment{keywords}{\noindent\small {\it Keywords\/}:}{\vskip 4pt}
\newenvironment{classification}{\noindent\small 2000 {\it Mathematics Subject
Classification\/}:}{\vskip 12pt}

\newcommand{\free}{{\mathbb F}}

\newcommand{\tensor}{\otimes}

\newcommand{\Tensor}{\hat{\otimes}}

\newcommand{\cstar}{{C^\ast}}

\newcommand{\A}{{\mathfrak A}}

\theoremstyle{plain}
\newtheorem*{lemma}{Lemma}
\newtheorem*{theorem}{Theorem}
\newtheorem*{corollary}{Corollary}
\theoremstyle{definition}
\newtheorem{definition}{Definition}
\theoremstyle{remark}
\newtheorem*{remark}{Remark}
\title{Biflatness and biprojectivity of the Fourier algebra}
\author{\textit{Volker Runde}}
\date{}
\begin{document}
\maketitle
\begin{abstract}
We show that the biflatness---in the sense of A.\ Ya.\ Helemski\u{\i}---of the Fourier algebra $A(G)$ of a locally compact group $G$ forces $G$ to either have an abelian subgroup of finite index or to be non-amenable without containing $\free_2$ as a closed subgroup. An analogous dichotomy is obtained for biprojectivity.
\end{abstract}

\section*{Introduction}
Abstract harmonic analysis is the study of locally compact groups $G$ and the various Banach algebras associated with them, such as the group algebra $L^1(G)$. In \cite{Joh1}, B.\ E.\ Johnson defined the class of amenable Banach algebras to consist of those Banach algebras that satisfy a certain cohomological triviality condition and showed that $L^1(G)$ is an amenable Banach algebra if and only if $G$ is an amenable locally compact group in the usual sense (\cite{Pat}). At about the same time, A.\ Ya.\ Helemski\u{\i} started to systematically develop a homological algebra with functional analytic overtones (see \cite{Hel} for an account). Amenability in the sense of \cite{Joh1} fits nicely into this framework and is closely related to the notion of biflatness. Another central notion in Helemki\u{\i}'s theory is that of biprojectivity. Like amenability, the biprojectivity of $L^1(G)$ singles out a natural class of groups: $L^1(G)$ is biprojective if and only if $G$ is compact (\cite[Theorem IV.5.13]{Hel}).
\par
The Fourier algebra $A(G)$---for arbitrary, not necessarily abelian $G$---was introduced by P.\ Eymard in \cite{Eym}. For abelian $G$, we have the dual group $\hat{G}$; in this case, $A(G)$ is nothing but $L^1(\hat{G})$ via the Fourier transform. For non-abelian $G$, however, $A(G)$ often displays a behavior strikingly differently from $L^1(G)$. For instance, $L^1(G \times H) \cong L^1(G) \tensor^\gamma L^1(H)$, with $\tensor^\gamma$ denoting the projective tensor product of Banach spaces, holds isometrically isomorphically for all $G$ and $H$ whereas $A(G \times H) \cong A(G) \tensor^\gamma A(H)$ holds isomorphically only if $G$ or $H$ is almost abelian, i.e., has an abelian subgroup of finite index (\cite{Los}).
\par 
Shortly after the publication of \cite{Eym}, H.\ Leptin provided a characterization of amenability in terms of the Fourier algebra: $G$ is amenable if and only if $A(G)$ has a bounded approximate identity (\cite{Lep}). Since amenable Banach algebras have bounded approximate identities, this suggests that $A(G)$---just like $L^1(G)$---is amenable if and only if $G$ is amenable, a view that seemed to have been widespread among mathematicians until the early 1990s. Then, in \cite{Joh2}, Johnson showed that there are compact groups $G$, such as $\mathrm{SO}(3)$, for which $A(G)$ fails to be amenable. On the positive side, it is not too difficult to see that $A(G)$ is indeed amenable if $G$ is almost abelian (\cite[Theorem 4.1]{LLW}). Eventually, B.\ E.\ Forrest and the author showed that $G$ being almost abelian is not only sufficient but also necessary for $A(G)$ to be amenable (\cite[Theorem 2.3]{FR}).
\par 
Since $A(G)$ is the predual of the von Neumann algebra generated by the left regular representation of $G$, it is an operator space in a canonical manner (see \cite{ER} for the theory of operator spaces). As it turns out, $A(G)$ is way better behaved as a completely contractive Banach algebra than as a mere Banach algebra. For instance, $A(G \times H) \cong A(G) \Tensor A(H)$ holds completely isometrically isomorphically for all $G$ and $H$, where $\Tensor$ is the projective tensor product of operator spaces, by \cite[Theorem 7.2.4]{ER}. Johnson's definition of an amenable Banach algebra can be modified to take operator space structures into account, which leads to the notion of operator amenability: this was done in \cite{Rua}, where Z.-J.\ Ruan showed that $A(G)$ is operator amenable if and only if $G$ is amenable. More generally, Helemski\u{\i}'s Banach homology can be developed in an operator space context as well (see \cite{Ari}, for instance), which leads to further interesting results about $A(G)$: it is operator biprojective if and only if $G$ is discrete (\cite{Ari}, \cite{RX}, or \cite{Woo}), and it is operator biflat for all so-called $\mathrm{[SIN]}$-groups (\cite{RX}) and possibly for all $G$ (\cite{ARS}).
\par 
In this brief note, we take a look at $A(G)$ in the framework of Helemski\u{\i}'s original Banach homology. We are interested in the properties of $G$ that are implied by the biflatness and biprojectivity of $A(G)$, respectively. Of course, $A(G)$ is biflat if $G$ is almost abelian. Our main result is that if $A(G)$ is biflat and $G$ is \emph{not} almost abelian, then $G$ is non-amenable, but does not contain $\free_2$, the free group in two generators, as a closed subgroup.
\section*{The result}
We begin with recalling the definitions of biprojectivity, biflatness, and amenability.
\par 
Given a Banach algebra $\A$, we let $m \!: \A \tensor^\gamma \A \to \A$ denote the multiplication map, i.e., $m(a \tensor b) =ab$ for $a,b \in \A$; if we want to emphasize the algebra, we sometimes write $m_\A$. The tensor product $\A \tensor^\gamma \A$ becomes a Banach $\A$-bimodule via
\[
  a \cdot (x \tensor y) := ax \tensor y \quad\text{and}\quad (x \tensor y) \cdot a := x \tensor ya \qquad (a,x,y \in \A)
\]
turning $m$ into a bimodule homomorphism.
\begin{definition} \label{homdefs}
Let $\A$ be a Banach algebra. Then $\A$ is called
\begin{enumerate}[\rm (a)]
\item \emph{biprojective} if $m \!: \A \tensor^\gamma \A \to \A$ has a bounded bimodule right inverse,
\item \emph{biflat} if $m^\ast \!: \A^\ast \to (\A \tensor^\gamma \A)^\ast$ has a bounded bimodule left inverse, and
\item \emph{amenable} if it is biflat and has a bounded approximate identity.
\end{enumerate}
\end{definition}
\begin{remark}
These definitions are not the original ones due to Johnson and Helemski\u{\i}, respectively, but are equivalent to them (\cite{Hel}).
\end{remark}
\par 
In order to deduce consequences for the structure of $G$ from the biflatness or biprojectivity of $G$, respectively, we require the following definition due to M.\ Leinert (\cite{Lei}):
\begin{definition} \label{Leidef}
Let $G$ be a discrete group. A subset $E$ of $G$ is called a \emph{Leinert set} if $\chi_E A(G) \cong \ell^2(E)$ holds isomorphically.
\end{definition}
\par 
Here, $\chi_E$ denotes the indicator function of $E$.
\par 
Trivially, all finite sets are Leinert sets. Remarkably, however, $\free_2$ contains \emph{infinite} Leinert set: this is the main result of \cite{Lei}. As was observed in \cite{Lei}, this implies that $A(\free_2)$ does not factor, i.e., there are functions in $A(\free_2)$ that are not product of two functions in $A(\free_2)$. More is true:
\begin{lemma}
Let $G$ be a locally compact group containing $\free_2$ as a closed subgroup. Then $m \!: A(G) \tensor^\gamma A(G) \to A(G)$ is not surjective.
\end{lemma}
\begin{proof}
Assume that $m_{A(G)}$ is surjective.
\par 
For any closed $F \subset G$, let $\rho_F$ denote the restriction map from $A(G)$, i.e., a function in $A(G)$ is restricted to $F$. By \cite[Theorem 16]{Her}, $\rho_{\free_2}$ maps $A(G)$ onto $A(\free_2)$. Since the diagram
\[
  \begin{diagram}
  A(G) \tensor^\gamma A(G)                     & \rTo^{m_{A(G)}}       & A(G) \\
  \dTo^{\rho_{\free_2} \tensor \rho_{\free_2}} &                       & \dTo_{\rho_{\free_2}} \\
  A(\free_2) \tensor^\gamma A(\free_2)         & \rTo^{m_{A(\free_2)}} & A(\free_2)
  \end{diagram}
\]
commutes, and since its columns are surjective, we conclude that $m_{A(\free_2)}$ is surjective as well.
\par 
Let $E \subset \free_2$ be an infinite Leinert set. By Definition \ref{Leidef}, is is clear that $\rho_E$ maps $A(\free_2)$ onto $\ell^2(E)$. Again, we have a commutative diagram
\[
  \begin{diagram}
  A(\free_2) \tensor^\gamma A(\free_2) & \rTo^{m_{A(\free_2)}} & A(\free_2) \\
  \dTo^{\rho_E \tensor \rho_E}         &                       & \dTo_{\rho_E} \\
  \ell^2(E) \tensor^\gamma \ell^2(E)   & \rTo^{m_{\ell^2(E)}}  & \ell^2(E)
  \end{diagram}
\]
with surjective columns, so that $m_{\ell^2(E)}$ has to be surjective, too. Since the range of $m_{\ell^2(E)}$ is contained in $\ell^1(E) \subsetneq \ell^2(E)$ by the Cauchy--Schwarz inequality, this yields a contradiction.
\end{proof}
\begin{remark}
For $G = \free_2$, the Lemma was already obtained, but never published, by H.\ Steiniger (\cite{Stei}) with a somewhat more technical proof that does not explicitly use Leinert sets. 
\end{remark}
\par 
Proving our main result now requires little more than assembling the right bits and pieces:
\begin{theorem}
Let $G$ be a locally compact group such that $A(G)$ is biflat. Then one of the following holds:
\begin{enumerate}[\rm (a)]
\item $G$ is almost abelian.
\item $G$ does not contain $\free_2$ as a closed subgroup, but fails to be amenable.
\end{enumerate} 
\end{theorem}
\begin{proof}
Suppose that $G$ is amenable. Then $A(G)$ has a bounded approximate identity by \cite{Lep}, making $A(G)$ amenable. By \cite[Theorem 2.3]{FR}, this means that (a) holds.
\par 
Suppose that $G$ is not amenable, and assume towards a contradiction that $G$ contains $\free_2$ as a closed subgroup. By the definition of biflatness, $m^\ast \!: A(G)^\ast \to (A(G) \tensor^\gamma A(G))^\ast$ has a bounded left inverse and thus, in particular, is injective with closed range. Consequently, $m \!: A(G) \tensor^\gamma A(G) \to A(G)$ has to be surjective, which contradicts the Lemma. Hence, (b) must hold.
\end{proof}
\begin{remark} 
The question of whether or not (discrete) groups as in (b) exist was open for several decades, and the belief that no such groups exist became known as ``von Neumann's conjecture''. Eventually, A.\ Yu.\ Ol'shanski\u{\i} came up with a counterexample to this conjecture (\cite{Ols}), so that condition (b) is not vacuous.
\end{remark}
\par 
As biprojectivity is stronger than biflatness, the conclusions of the Theorem apply, in particular, if $A(G)$ is biprojective. Furtheremore, by general Banach algebra theory (\cite[Corollary 2.8.42]{Dal}), the biprojectivity of $A(G)$ forces $G$ to be discrete. 
\par
We summarize:
\begin{corollary}
Let $G$ be a locally compact group group such that $A(G)$ is biprojective. Then $G$ is discrete, and one of the following holds:
\begin{enumerate}[\rm (a)]
\item $G$ is almost abelian.
\item $G$ does not contain $\free_2$ as a closed subgroup, but fails to be amenable.
\end{enumerate} 
Furthermore, for any discrete $G$, \emph{(a)} is also sufficient for $A(G)$ to be biprojective.
\end{corollary}
\begin{proof}
Only the ``furthermore'' part still needs consideration. Suppose that $G$ is discrete and almost abelian. Then $A(G)$ is amenable, and since $A(G)$ is Tauberian, the discreteness of $G$ forces multiplication in $A(G)$ to be compact. By \cite[Corollary 3.2]{LLRW}, this means that $A(G)$ is biprojective.
\end{proof}
\renewcommand{\baselinestretch}{1.0}
\renewcommand{\baselinestretch}{1.2}
\dated
\vfill
\begin{tabbing}
\textit{Author's address}: \= Department of Mathematical and Statistical Sciences \\
                           \> University of Alberta \\
                           \> Edmonton, Alberta \\
                           \> Canada T6G 2G1 \\[\medskipamount]
\textit{E-mail}:           \> \texttt{vrunde@ualberta.ca}\\[\medskipamount]
\textit{URL}:              \> \texttt{http://www.math.ualberta.ca/$^\sim$runde/}
\end{tabbing}
\end{document}